\documentclass[11pt]{amsart} 
\usepackage{url}
\newtheorem{theorem}{Theorem}[section]
\newtheorem*{sharkovskystheorem}{Sharkovsky's Theorem}
\begin{document} 
\noindent 
{\small Topology Atlas Invited Contributions {\bf 6} no.~1 (2001) 4 pp.} 
\vspace{\baselineskip}
\title{Combinatorial Dynamics} 
\author{Micha\l\ Misiurewicz} 
\address{Department of Mathematical Sciences,
Indiana University--Purdue University Indianapolis, 
402 N.~Blackford Street, Indianapolis, IN 46202-3216, USA}
\email{mmisiure@math.iupui.edu}
\urladdr{http://www.math.iupui.edu/$\tilde{~}$mmisiure/}
\thanks{Note: Detailed treatment of Combinatorial Dynamics in dimension 
$1$ (with an extensive bibliography and historical remarks) can be found
in the book \cite{1}. For dimension $2$, see for instance
\cite{2,3,4,5}.}
\maketitle

\section{Interval}

Combinatorial Dynamics has its roots in Sharkovsky's Theorem. 
This beautiful theorem describes the possible sets of periods of all 
cycles of a continuous map of an interval (or the real line) into itself.
Here by a {\em cycle} I mean a periodic orbit, and by a
{\em period} its minimal period.

Consider the following {\em Sharkovsky's ordering}
$<_s$ of the set $\mathbb{N}$ of natural numbers:
\begin{align*}
1 <_s 2 <_s 2^2 <_s 2^3&<_s 2^4 <_s \cdots <_s 2^n <_s \cdots\\
&\vdots
\end{align*}
$$
\begin{array}{lllll}
\dots <_s 11 \times 2^n&
<_s 9 \times 2^n&
<_s 7 \times 2^n&
<_s 5 \times 2^n&
<_s 3 \times 2^n\\
&
&
\vdots&
&
\\
\dots <_s 11 \times 2^2&
<_s 9 \times 2^2&
<_s 7 \times 2^2&
<_s 5 \times 2^2&
<_s 3 \times 2^2\\
\dots <_s 11 \times 2&
<_s 9 \times 2&
<_s 7 \times 2&
<_s 5 \times 2&
<_s 3 \times 2\\
\dots <_s 11&   
<_s 9&
<_s 7&
<_s 5&
<_s 3 
\end{array}
$$

Denote by $S(n)$ the {\em initial segment} of this ordering ending at 
$n$ in $\mathbb{N}$, that is 
$S(n) = \{m \in \mathbb{N} : m <_s n \} \cup \{n\}$. 
We also set $S(2^\infty) = \{1,2,2^2,2^3,2^4,\dots\}$. 
For a map $f\colon X \to X$ we denote by $P(f)$ the set of 
periods of all cycles of $f$.

\begin{sharkovskystheorem}
Let $I$ be a closed interval.
Then for every continuous map $f\colon I \to I$ there exists $n$ in 
$\mathbb{N} \cup \{2^\infty\}$ such that $P(f) = S(n)$. Conversely, for 
every $n$ in $\mathbb{N} \cup \{2^\infty\}$ there exists a continuous map
$f\colon I \to I$ such that $P(f) = S(n)$.
\end{sharkovskystheorem}

The same holds for continuous maps of the real line into itself,
except that we may have $P(f) = \emptyset$.

All proofs of Sharkovsky's Theorem inevitably lead to an idea of a 
``type'' of a cycle. For instance, when ordering the points of a cycle
$p_1 < p_2 < \cdots < p_n$, one gets a cyclic permutation $\sigma$ 
corresponding to it, such that $p_i$ is mapped to 
$p_{\sigma(i)}$.

Combinatorial Dynamics deals with types of cycles (defined in various ways 
and not only for interval maps), their coexistence, complexity (usually 
measured by the topological entropy of the system) and properties of maps 
with cycles of a given type. The word ``Combinatorial'' refers not only to 
the fact that we often deal with permutations, but also to the techniques 
involved in the proofs, that are often of a combinatorial nature.

Another way of looking at the permutation given by a cycle (and forgetting
about the orientation of the interval) is to say that two cycles, $P$ of
the map $f\colon I \to I$ and $Q$ of the map $g\colon J \to J$, have the 
same {\em pattern} if there is a homeomorphism $h\colon I \to J$ such that
$h(P) = Q$ and $h(f(x)) = g(h(x))$ for $x$ in $P$. We say that a pattern
$A$ {\em forces} pattern $B$ if any interval map with a cycle of pattern
$A$ has also a cycle of pattern $B$. For other definitions of the ``type''
of a cycle, we define forcing in the same way. Note that by Sharkovsky's
Theorem, forcing among periods is a linear ordering.

\begin{theorem}
Forcing among patterns is a partial ordering.
\end{theorem}

There is another notion, between the period and the pattern, that gives us
also a linear ordering. Let $q$ be the period of a cycle $P$ of $f$, and
let $m$ be the number of times $f(x) - x$ changes sign when $x$ moves
along $P$ (that is, $x = y, f(y), \dots, f^{q-1}(y), y$). Note that
$p = m/2$ is an integer. We call the pair $(p,q)$ the {\em over-rotation
pair} of $P$ and the number $p/q$ the {\em over-rotation number} of $P$.

\begin{theorem}
For a given continuous interval map, the set of over-rotation numbers of 
all cycles of period larger than $1$ is either the intersection of 
$\mathbb{Q}$ with an interval whose one endpoint is $1/2$, or $\{1/2\}$, 
or $\emptyset$.
\end{theorem}

\begin{theorem}
Forcing among over-rotation pairs is a linear ordering.
\end{theorem}

This ordering can be described as follows. Suppose we want to compare two
over-rotation pairs $(p,q)$ and $(r,s)$. If $p/q > r/s$ then $(p,q)$
forces $(r,s)$; if $p/q < r/s$ then $(r,s)$ forces $(p,q)$. If
$p/q = r/s$ then we write $p/q$ as a fraction in the reduced form and take
its denominator $k$. Both $q,s$ are divisible by $k$, and the relation
between $(p,q)$ and $(r,s)$ is the same as between $q/k$ and $s/k$ in the
Sharkovsky's ordering.

Often we can identify the simplest patterns within some classes, for
instance patterns that do not force any other pattern of the same period
({\em primary patterns}) or patterns that do not force any other pattern
of the same over-rotation number ({\em twist patterns}).

The infimum of the topological entropies of maps with a cycle of a given
pattern is called the entropy of a pattern. It indicates how chaotic a map
has to be to accommodate a cycle of a given pattern. There are simple ways
to compute this entropy, and cycles of the smallest entropy with a given
period or over-rotation number have been identified.

\section{Other one-dimensional spaces}

Naturally, attempts have been made to built similar theories for maps of
other one-dimensional spaces. The first obvious candidate is the circle.
Here continuous maps are immediately classified by their degree, and it
turns out that the most interesting case is when the degree is $1$ (maps
homotopic to the identity). In this case, an important characteristics of
a cycle is its {\em rotation number}. It is defined as the number of times
the orbit goes around the circle, divided by the period. To measure it, we
have to fix a lifting $F\colon \mathbb{R} \to \mathbb{R}$ of the map,
where the circle is $\mathbb{R} / \mathbb{Z}$, and then it is equal to the
average along the cycle of the displacement function $\varphi(x) = f(y) - y$
where $y$ in $\mathbb{R}$ projects to $x$. The over-rotation number for
interval maps mentioned in the preceding section is an analogue of this
notion; instead of the displacement function $\varphi$ we took $\psi(x)$ 
equal to $1/2$ if $(f^2(x) - f(x))(f(x) - x) > 0$ and $0$ otherwise.

\begin{theorem}
For a given continuous circle map of degree $1$, the set of rotation 
numbers of all cycles is either the intersection of $\mathbb{Q}$ with a 
closed interval, or a singleton, or $\emptyset$.
\end{theorem}

Possible sets of periods are obtained by choosing the endpoints of the
rotation interval and for each rational endpoint one element of
$\mathbb{N} \cup \{2^\infty\}$. Then periods are all denominators of the
fractions from the interior of the rotation interval (not necessarily in
the reduced form) and initial segment(s) of the Sharkovsky's ordering
times the denominators of the endpoints of the rotation interval
(reduced).

Thus, the description of the set of periods of such a map involves more
data than in the case of interval maps, so the patterns and the forcing
relation do not play important role here. Minimal entropy of maps with a
given rotation interval is known.

For continuous tree maps (except maps of $n$-ods which fix the central
point, where the situation resembles that for interval maps), despite
substantial progress, the situation is still unclear. It seems that the
idea of patterns can work, but one cannot fix the tree. That is, maps of
different trees can have cycles with the same pattern.

For graph maps, basically nothing is known from the point of view of
Combinatorial Dynamics.

\section{Dimension 2}

Usually notions that make sense for noninvertible maps in dimension $1$
have their counterparts for homeomorphisms in dimension $2$. For
noninvertible maps in dimension $2$ or homeomorphisms in dimension $3$,
there are already too much ``degrees of freedom'', allowing us to
construct various counterexamples.

For a map of a disk $D$ into itself, one can define a pattern of a cycle
in the following way. Cycles $P$ of $f$ and $Q$ of $g$ have the same
pattern if there exists a homeomorphism $h \colon D \to D$ such that
$h(P) = Q$, $h(f(x)) = g(h(x))$ for $x$ in $P$ and $g$ is isotopic to $h
\circ f \circ h^{-1}$ modulo $Q$ (that is, the isotopy maps act on the
points of $Q$ in the same way as $g$). This is equivalent to the condition
that the orbits that are obtained from $P$ and $Q$ when we take
suspensions of $f$ and $g$ respectively have the same braid types (modulo
full twists). The same definitions work for maps of any surface map
isotopic to the identity.

Forcing among patterns can be defined in the same way as for interval
maps, and Theorem 1.1 holds also in this case. Similar questions can be
asked, although usually they are much harder than for interval maps. The
classification of patterns uses the Nielsen-Thurston theory.

While for interval maps it is obvious whether two given cycles have the
same pattern, in two dimensions the algorithm for doing that is quite
complicated.

Rotation theory has a nice generalization for homeomorphisms of the torus
isotopic to the identity. The displacement function takes values in
$\mathbb{R}^2$, so this is where the rotation vectors live.

\begin{theorem}
For a homeomorphism of the torus isotopic to the identity the set of 
rotation vectors of cycles is equal to the intersection of $\mathbb{Q}^2$
with some compact convex set, except perhaps some points from interiors 
of segments contained in the boundary of this set.
\end{theorem}

An interesting problem is whether any compact convex subset of the
plane is the rotation set of a torus homeomorphism.

\end{document}